# Tropical Hurwitz numbers

**Renzo Cavalieri · Paul Johnson · Hannah Markwig**

Received: 24 February 2009 / Accepted: 4 December 2009 / Published online: 24 December 2009

**Abstract** *Hurwitz numbers* count genus $g$, degree $d$ covers of $\mathbb{P}^1$ with fixed branch locus. This equals the degree of a natural *branch map* defined on the Hurwitz space. In *tropical geometry*, algebraic curves are replaced by certain piece-wise linear objects called tropical curves. This paper develops a tropical counterpart of the branch map and shows that its degree recovers classical Hurwitz numbers. Further, the combinatorial techniques developed are applied to recover results of Goulden et al. (in Adv. Math. 198:43–92, 2005) and Shadrin et al. (in Adv. Math. 217(1):79–96, 2008) on the piecewise polynomial structure of double Hurwitz numbers in genus 0.

P. Johnson was supported in part by University of Michigan RTG grant 0602191.

H. Markwig was supported by the German Research Foundation (Deutsche Forschungsgemeinschaft (DFG)) through the Institutional Strategy of the University of Göttingen.

R. Cavalieri (✉)
Department of Mathematics, Colorado State University, Weber Building, Fort Collins, CO 80523-1874, USA
e-mail: renzo@math.colostate.edu

P. Johnson
Department of Mathematics, University of Michigan, 2074 East Hall, 530 Church Street, Ann Arbor, MI 48109-1043, USA
e-mail: pdjohnso@umich.edu

H. Markwig
CRC “Higher Order Structures in Mathematics”, Georg August Universität Göttingen, Bunsenstr. 3-5, 37073 Göttingen, Germany
e-mail: hannah@uni-math.gwdg.de

## 1 Introduction

Hurwitz numbers are important objects connecting the geometry of algebraic curves to the combinatorics of the symmetric group. Geometrically, Hurwitz numbers count genus $g$, degree $d$ covers of $\mathbb{P}^1$, with specified ramification profile over a fixed set of $n$ points in $\mathbb{P}^1$. By matching a cover with a (equivalence class of) monodromy representation, such count is equivalent to choices of $n$-tuples of elements of $S_d$ multiplying to the identity element and acting transitively on the set $\{1, \dots, d\}$.

This connection dates back to Hurwitz himself, and has provided a rich interplay between the two fields for a long time. In more recent times, Hurwitz numbers have found a prominent role in the study of the moduli space of curves and in Gromov-Witten theory. The moduli space of genus $g$, degree $d$ covers of $\mathbb{P}^1$ with only simple ramification (Hurwitz space) admits a natural *branch map* recording the position of the branch points on $\mathbb{P}^1$. The degree of the branch map onto its image is tautologically equal to a Hurwitz number. The Hurwitz space sits inside the moduli space of stable maps $\overline{M}_g(\mathbb{P}^1, d)$. However, the latter is a singular, non-equidimensional stack. In [3], Fantechi and Pandharipande define a branch map

$$br : \overline{M}_g\big(\mathbb{P}^1, d\big) \to \mathbb{P}^n = \mathrm{Sym}^n\big(\mathbb{P}^1\big)$$

and show that the virtual degree

$$br^{-1}(pt) \cap \big[\overline{M}_g\big(\mathbb{P}^1, d\big)\big]^{\mathrm{vir}}$$

still recovers the appropriate Hurwitz number.

The theory of relative stable maps [9, 10] extends this scenario to more general Hurwitz numbers, with arbitrary ramification profiles over the branch points. Connections with the moduli spaces of curves, the most remarkable being the *ELSV formula* [2], are produced by virtual localization on these moduli spaces [7].

In tropical geometry, algebraic curves are degenerated to certain piece-wise linear graphs called tropical curves. This process "loses a lot of information", but many properties of the algebraic curve can be read off the tropical curve, and many theorems that hold for algebraic curves remarkably continue to hold on the tropical side. One of the fields in which tropical geometry has had significant success recently is enumerative geometry. Moduli spaces of curves and maps are important objects for the study of enumerative geometry, and so their tropical counterparts have been at the center of attention.

The aim of this paper is twofold: first, we develop tropical branch maps and show that their degrees are equal to the Hurwitz numbers, just as in the "classical world". Second, we use one of lemmas to understand the combinatorial structure of double Hurwitz numbers; in this paper we develop this application in genus 0.

As a first step, we understand degree $d$ covers of $\mathbb{P}^1$ tropically. It is natural to think of them as elements in tropical $M_{g,0}(\mathbb{P}^1, d)$. An element in tropical $M_{g,0}(\mathbb{P}^1, d)$ roughly consists of a graph of genus $g$ and a map to tropical $\mathbb{P}^1$ satisfying certain conditions (see Sect. 5.1). We think of tropical $\mathbb{P}^1$ as $\mathbb{R} \cup \{\pm\infty\}$ (see [11]). Next, we seek for a tropical counterpart of a branch point. Thinking of a branch point as a point where several sheets of the fibers of the map come together, it is natural to interpret a tropical branch point as a (more than 2-valent) vertex of the abstract tropical curve mapping to $\mathbb{P}^1$. So far, this notion seems to work effectively only in the generic case of simple ramification, corresponding to trivalent vertices.

Luckily, we have the freedom to impose arbitrary ramification profiles over the two special points $\pm\infty$. In these case the ramification data is encoded in the collection of weights of the edges "going to $\infty$".

To sum up, we think of a tropical cover of $\mathbb{P}^1$ as an element in $M_{g,0,\mathrm{trop}}(\mathbb{P}^1, d)$; the branch map assigns the positions of the vertices. To adapt the notion of tropical moduli spaces of maps to our situation, we introduce a labeling on the vertices. The moduli space we define is a weighted polyhedral complex. Although we would prefer to give it the structure of a tropical variety, the structure of polyhedral complex is enough to work with the branch map. We adopt this shortcut because for higher genus it is not yet known how the tropical moduli spaces can be understood as tropical varieties. The degree of the tropical branch map is well defined.

Our main theorem (Theorem 5.28 in Sect. 5.3) shows that the degree of the tropical branch map is equal to the classical double Hurwitz numbers.

The strategy of proof is natural. We interpret each graph on the tropical side as a family of monodromy representations contributing to the classical Hurwitz number. The number of elements in this family is determined via a cut and join analysis, and coincides with the *tropical multiplicity* of the graph: how many times a tropical curve with that graph occurs in the preimage of the branch map times the multiplicity of the branch map.

Lemma 4.2, which shows how cut and join recursion is conveniently organized as a sum over weighted graphs, leads to new and elementary proofs of existing results about the combinatorial structure of double Hurwitz numbers. It was shown in [6] that double Hurwitz numbers are piecewise polynomial in the entries of the two partitions, and in [13] this structure was investigated in genus 0. In particular, it was shown that the chambers of polynomiality are delimited by hyperplanes parameterizing Hurwitz data allowing the existence of disconnected cover curves, and that the way the Hurwitz polynomials vary across a wall can be expressed recursively in terms of the Hurwitz data for the connected components of these disconnected covers. Using Lemma 4.2, we give an elementary proof of the fact that Hurwitz numbers are piecewise polynomial in genus 0, and in Theorem 6.10 we give a proof of the wall-crossing formula. Although these proofs do not logically depend on tropical geometry, Lemma 4.2 was only discovered with tropical geometry as a motivation.

The same approach to polynomiality and wall crossing can be extended to higher genus, where wall-crossing formulas were previously unknown. However, the combinatorics required are considerably more sophisticated. These results are presented in [1].

The paper is organized as follows. We first recall the definition of Hurwitz numbers (Sect. 2) and discuss the cut and join equations (Sect. 3). In Sect. 4, we deduce a weighted graph count which computes the Hurwitz number. Section 5 establishes tropical Hurwitz theory. In Sect. 5.1 we define the tropical moduli space and in Sect. 5.2 the tropical branch map. Our main theorem that the degree of the tropical branch map equals the Hurwitz number is formulated and proved in Sect. 5.3. Section 6 explores combinatorial properties of double Hurwitz numbers in genus 0. In Sect. 6.1 we show piecewise polynomiality and identify the polynomiality chambers. In Sect. 6.2 we prove the wall-crossing formula.

## 2 Hurwitz numbers

In this section we recall the definition and some basic facts about Hurwitz numbers.

**Definition 2.1** Fix $r+s$ points $p_1,\dots,p_r,q_1,\dots,q_s$ on $\mathbb{P}^1$, and $\eta_1,\dots,\eta_r$ partitions of the integer $d$. The *Hurwitz number*:

$$H_d^g(\eta_1,\dots,\eta_r)$$

$$:= \text{weighted number of} \left\{ \begin{array}{l} \text{degree } d \text{ covers } C \xrightarrow{\pi} \mathbb{P}^1 \text{ such that:} \\ \bullet\ C \text{ is a smooth connected curve of genus } g; \\ \bullet\ \pi \text{ is unramified over } \mathbb{P}^1 \setminus \{p_1,\dots,p_r,q_1,\dots,q_s\}; \\ \bullet\ \pi \text{ ramifies with profile } \eta_i \text{ over } p_i; \\ \bullet\ \pi \text{ has simple ramification over } q_i. \end{array} \right\}$$

Each cover $\pi$ is weighted by $1/|\operatorname{Aut}(\pi)|$.

Note that this is independent of the locations of the $p_i$ and $q_i$. For a partition $\eta$, let $\ell(\eta)$ denote the number of parts of $\eta$. By the Riemann–Hurwitz formula, we have that

$$2-2g = 2d - dr - s + \sum_{i=1}^{r} \ell(\eta_i),$$

and hence $s$ is determined by $g, d$ and $\eta_1,\dots,\eta_r$.

It is often common language to use *Hurwitz number* for the generic case $H_d^g$ when all ramification is simple; *simple* (resp. *double*) Hurwitz number when one (resp. two) point of arbitrary ramification are prescribed.

A ramified cover is essentially equivalent information to a monodromy representation it induces; thus, an equivalent definition of Hurwitz number counts the number of homomorphisms $\varphi$ from the fundamental group $\Pi_1$ of $\mathbb{P}^1 \setminus \{p_1,\dots,p_r,q_1,\dots,q_s\}$ to the symmetric group $S_d$ such that:

- the image of a loop around $p_i$ has cycle type $\eta_i$;
- the image of a loop around $q_i$ is a transposition;
- the subgroup $\varphi(\Pi_1)$ acts transitively on the set $\{1,\dots,d\}$.

This number is divided by $|S_d|$, to account both for automorphisms and for different monodromy representations corresponding to the same cover.

## 3 Cut and join

The *Cut and Join equations* are a collection of recursions among Hurwitz numbers. In the most elegant and powerful formulation they are expressed as one differential operator acting on an appropriate potential function. Since our use of cut and join is unsophisticated, we limit ourselves to a basic discussion, and refer the reader to [5] for a more in-depth presentation.

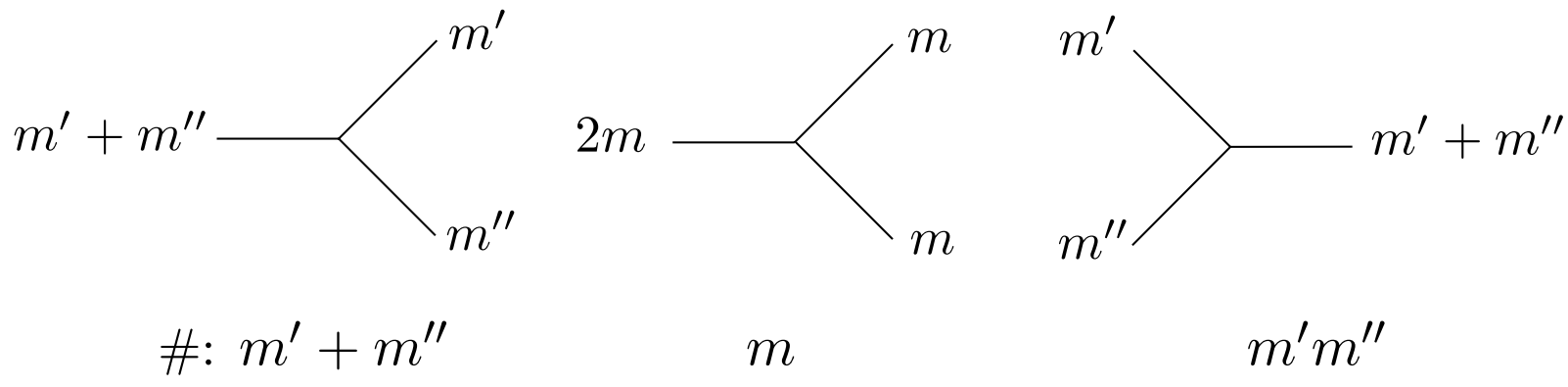

**Fig. 1** Composing with a transposition in $S_d$. How it affects the cycle type of $\sigma$ and multiplicity

Let $\sigma \in S_d$ be a fixed element of cycle type $\eta = (n_1, \dots, n_l)$, written as a composition of disjoint cycles as $\sigma = c_l \cdots c_1$. Let $\tau = (ij) \in S_d$ vary among all transpositions. The cycle types of the composite elements $\tau\sigma$ are described below.

- **cut:** if $i, j$ belong to the same cycle (say $c_l$), then this cycle gets "cut in two": $\tau\sigma$ has cycle type $\eta' = (n_1, \dots, n_{l-1}, m', m'')$, with $m' + m'' = n_l$. If $m' \neq m''$, there are $n_l$ transpositions giving rise to an element of cycle type $\eta'$. If $m' = m'' = n_l/2$, then there are $n_l/2$.
- **join:** if $i, j$ belong to different cycles (say $c_{l-1}$ and $c_l$ ), then these cycles are "joined": $\tau\sigma$ has cycle type $\eta' = (n_1, \dots, n_{l-1} + n_l)$. There are $n_{l-1}n_l$ transpositions giving rise to cycle type $\eta'$.

*Example 3.1* Let $d = 4$. There are 6 transpositions in $S_4$. If $\sigma = (12)(34)$ is of cycle type $(2, 2)$, then there are 2 transpositions ((12) and (34)) that "cut" $\sigma$ to give rise to a transposition and $2 \cdot 2$ transpositions ((13), (14), (23), (24)) that "join" $\sigma$ into a four-cycle.

For readers allergic to notation, Fig. 1 illustrates the above discussion.

## 4 Double Hurwitz numbers and Weighted graph sums

The analysis in Sect. 3 leads us to compute double Hurwitz numbers in terms of a weighted sum over graphs. The idea is to start at one of the special points, and count all possible monodromy representations as each transposition gets added until one gets to the second special point with the specified cycle type. We now make this precise.

Fix $g$ and let $\eta = (n_1, \dots, n_k)$ and $\nu = (m_1, \dots, m_l)$ be two partitions of $d$. Denote by $s = 2g - 2 + l + k$ the number of non-special branch points, determined by the Riemann–Hurwitz formula.

**Definition 4.1** *Monodromy graphs* project to the segment $[0, s + 1]$ and are constructed according to the following procedure:

(a) Start with $k$ small segments over 0 labeled $n_1, \dots, n_k$. We call these $n$'s the weights of the strands.

(b) Over the point 1 create a three-valent vertex by either joining two strands or splitting one with weight strictly greater than 1. In case of a join, label the new strand with the sum of the weights of the edges joined. In case of a cut, label the two new strands in all possible (positive) ways adding to the weight of the split edge.
(c) Consider only one representative for any isomorphism class of labeled graphs.
(d) Repeat (b) and (c) for all successive integers up to $s$.
(e) Retain all connected graphs that "terminate" with $l$ points of weight $m_1, \dots, m_l$ over $s+1$.

*Note* The graphs constructed above should be considered as abstract graphs with weighted edges and a map to the segment $[0, s+1]$. In other words, the relative positions of the strands is irrelevant, and there are no crossings between the strands.

**Lemma 4.2** *The double Hurwitz number $H^g_d(\eta, \nu)$ is computed as a weighted sum over monodromy graphs. Each monodromy graph is weighted by the product of the following factors*:

(i) *The number $\epsilon(\eta)$ of elements of $S_d$ of cycle type $\eta$.*
(ii) $|\mathrm{Aut}(\eta)|$.
(iii) *For every vertex*, *the product of the degrees of edges coming into the vertex from the left.*
(iv) *A factor of* $1/2$ *for any balanced fork or wiener.*
(v) $1/d!$.

A *balanced fork* is a tripod with weights $n, n, 2n$ such that the vertices of weight $n$ lie over 0 or $s+1$. A *wiener* consists of a strand of weight $2n$ splitting into two strands of weight $n$ and then re-joining. See Fig. 2.

*Proof* Recall how a ramified cover gives rise to a monodromy representation. Pick a point $x$ outside the branch locus in $\mathbb{P}^1$, and label the preimages $1, \dots, d$. Choose a set of loops based at $x$ winding around each branch point, letting the first and the last loop wind around the two special points. The liftings of the loops give rise to permutations of the preimages of $x$. The two special points give permutations $\sigma_\eta$, $\sigma_\nu$ of cycle type $\eta$ and $\nu$; all other points give transpositions $\tau_i$. The product

$$\sigma_\nu \tau_s \cdots \tau_1 \sigma_\eta$$

is the identity.

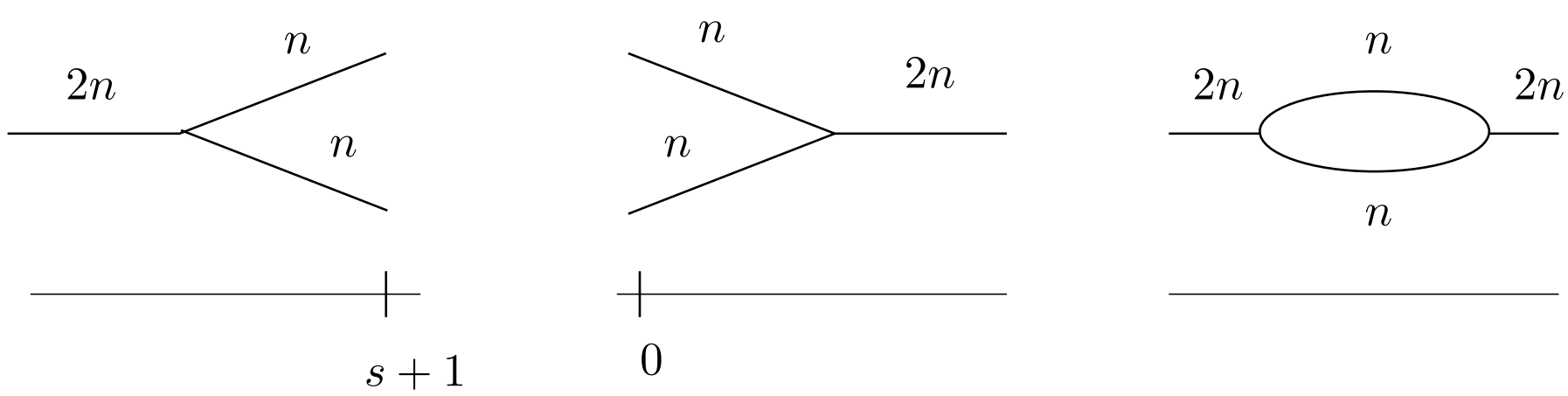

**Fig. 2** Balanced right pointing forks, balanced left pointing forks and wiener

The first permutation $\sigma_\eta$ can be chosen in $\epsilon(\eta)$ ways (i). We represent this cycle type as described in (a). The strands should be associated to the cycles of $\sigma_\eta$, but we are only labeling the strands by their weights: there are $|\mathrm{Aut}(\eta)|$ distinct ways of assigning a cycle of $\eta$ to each strand (ii). If two strands with the same multiplicity get joined to form a left pointing fork, then the two "teeth" of the fork are indistinguishable—this gives us a $1/2$ factor for every balanced left pointing fork (iv).

The analysis of Sect. 3 tells us how the cycle type changes every time a transposition is composed. The cut and join action is represented as described in (b), (c), (d); lastly, (e) ensures the connectedness of the cover and the right cycle type over the last point. The weights (iii) are precisely the multiplicities given by the cut and join analysis, with one exception. We count the cutting of $2n$ into $n, n$ twice as much as the cut and join equation prescribes. If the two $n$-strands have distinguishable evolution after the splitting, then it matters which cycle has which evolution: each $(n, n)$ cycle counts for 2. The two $n$ strands do not have a distinguishable evolution only in the case of balanced right pointing forks and wieners: in this case we want the original cut and join count. This gives a factor of $(1/2)^{\#\text{b.r.forks}+\#\text{wieners}}$ correcting our convention (iv). Finally we divide by $d!$ to account for the action of $S_d$ by conjugation, corresponding to a relabeling of the $d$ preimages of $x$ (v). □

*Remark 4.3* The genus of all graphs in the graph sum is $g$. This is immediately seen by combining the computation of the Euler characteristic of the graphs with the Riemann–Hurwitz formula.

We organize our graph sum so as to make it transparent that it counts monodromy representations. In our presentation the symmetry between the two special partitions is not obvious. However it suffices to notice that

$$\epsilon(\eta)\big|\mathrm{Aut}(\eta)\big| \cdot \prod_1^k n_i = d! \tag{1}$$

to recover the desired symmetry. With this substitution, we obtain the following.

**Corollary 4.4** *The formula for double Hurwitz numbers from Lemma 4.2 simplifies to*:

$$H_d^g(\eta, \nu) = \sum_{\Gamma} \frac{1}{|\mathrm{Aut}(\Gamma)|} \prod w(e), \tag{2}$$

*where we take the product of all the interior edge weights*; *the factors of* $1/2$ *coming from the balanced forks and wieners amount to the size of the automorphism group of our decorated graphs.*

*Example 4.5* We illustrate our procedure by computing:

$$H_4^1\big((4), (2, 2)\big) = 14.$$

Table 1 shows the type of contributing graphs and the various contributions discussed in this section, both in the form of Lemma 4.2 and of Corollary 4.4.

**Table 1**

| Graph type | (i) | (ii) | (iii) | (iv) | (v) | $\prod w(e)$ | $\lvert \mathrm{Aut}(\Gamma)\rvert$ | Total |
|---|---|---|---|---|---|---|---|---|
| 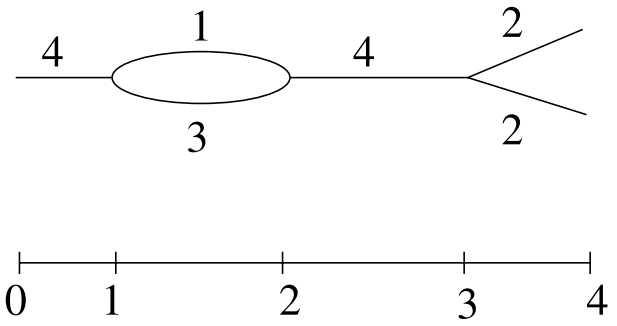 | 6 | 1 | 48 | $\frac{1}{2}$ | $\frac{1}{24}$ | 12 | 2 | 6 |
| 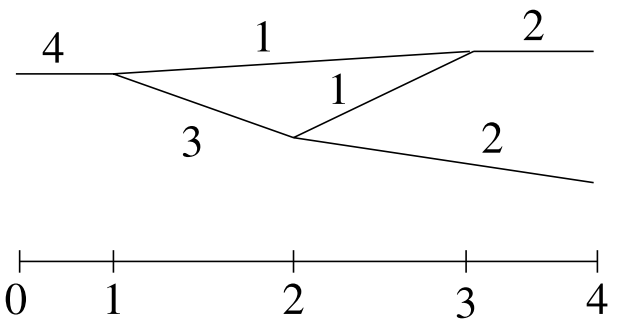 | 6 | 1 | 12 | 1 | $\frac{1}{24}$ | 3 | 1 | 3 |
| 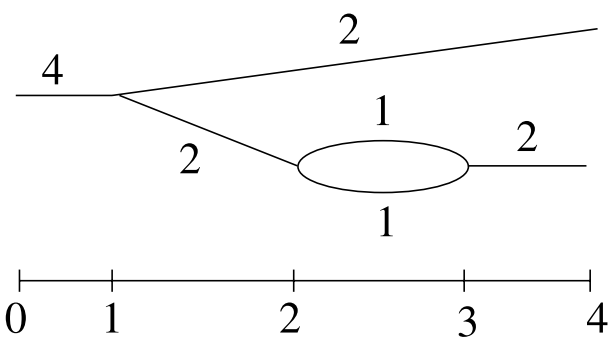 | 6 | 1 | 8 | $\frac{1}{2}$ | $\frac{1}{24}$ | 2 | 2 | 1 |
| 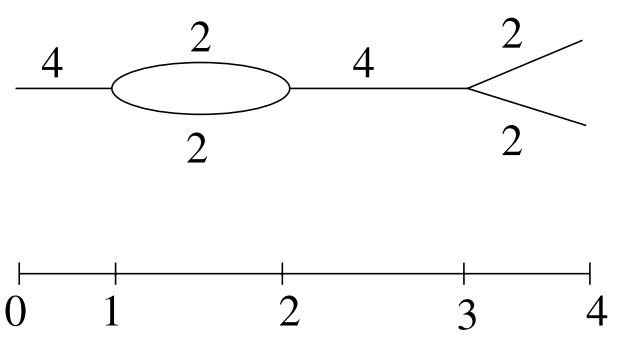 | 6 | 1 | 64 | $\frac{1}{4}$ | $\frac{1}{24}$ | 16 | 4 | 4 |

## 5 Tropical Hurwitz theory

### 5.1 Tropical maps to $\mathbb{P}^1$

In this section, we define the tropical moduli spaces needed for tropical Hurwitz numbers.

Let $\Gamma$ be a connected graph without 2-valent vertices. We call *ends* of $\Gamma$ edges adjacent to a 1-valent vertex. Edges which are not ends are called *bounded edges*. We denote the set of vertices by $\Gamma^0$, the subset of 1-valent vertices by $\Gamma^0_\infty$ and the subset of more than 1-valent vertices (called *inner vertices*) $\Gamma^0_0$. Likewise, the set of edges is $\Gamma^1$, the subset of ends $\Gamma^1_\infty$ and the bounded edges $\Gamma^1_0$. We call a pair $F = (V, e)$ where $e$ is an edge of $\Gamma$ and $V \in \partial e$ a *flag* of $\Gamma$ and think of it as a "directed edge", an edge pointing away from its end vertex $V$.

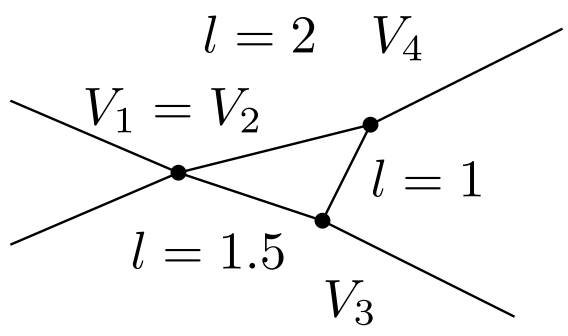

**Fig. 3** An abstract tropical curve

An *abstract tropical curve* is a connected graph $\Gamma$ without 2-valent vertices, whose edges $e$ are equipped with a length $l(e) \in \mathbb{R}_{>0} \cup \{\infty\}$. Ends $e \in \Gamma^1_\infty$ have length $l(e) = \infty$ and bounded edges $e \in \Gamma^1_0$ have a length $l(e) \in \mathbb{R}_{>0}$. We can think of the edges of an abstract tropical curves as intervals $(0, l(e))$.

The *genus* of an abstract tropical curve is the genus of $\Gamma$, equal to $h_1(\Gamma)$ since $\Gamma$ is connected.

The valency of a vertex is denoted by $\mathrm{val}(V)$. An *abstract tropical curve with g-labels* is an abstract tropical curve $\Gamma$ of genus $g' \leq g$ where each inner vertex $V$ is labeled with $\mathrm{val}(V) - 2 + 2k_V$ numbers for some $k_V \geq 0$, and points $p \in (0, l(e))$ on each edge $e$ are labeled with $2k_p$ numbers for some $k_p \geq 0$, such that the disjoint union of all labels equals $\{1, \ldots, s - 2 + 2g\}$ where $s$ is the number of ends.

For $i \in \{1, \ldots, s - 2 + 2g\}$ we denote by $V_i$ the vertex or point which has the label $i$.

Two abstract tropical curves with $g$-labels $\Gamma$ and $\tilde{\Gamma}$ are called isomorphic (and will from now on be identified) if there is a homeomorphism $\Gamma \to \tilde{\Gamma}$ mapping the label $i$ in $\Gamma$ to $i$ in $\tilde{\Gamma}$ for all $i$ and such that every edge of $\Gamma$ is mapped bijectively onto an edge of $\tilde{\Gamma}$ by an affine map of slope $\pm 1$, i.e. by a map of the form $t \mapsto a \pm t$ (where $a = 0$ or $a = l(e)$, and we again identify an edge of length $l(e)$ with the interval $(0, l(e))$).

The *combinatorial type* of an abstract tropical curve with $g$-labels is the data obtained when dropping the information about the lengths of the edges.

*Example 5.1* Figure 3 shows an abstract tropical curve of genus 1 with 1-labels. All ends have length $\infty$, the bounded edges have the length written next to them. The labels are at the inner vertices, marked with a black dot in the picture.

A graph of genus $g$ has $\#\Gamma^1_0 = \#\Gamma^1_\infty - 3 + 3g - \sum_{V \in \Gamma^0_0} (\mathrm{val}(V) - 3)$ bounded edges. We call a graph for which every vertex has valence 1 or 3 a 3-*valent graph*. In particular, a 3-valent graph has $\#\Gamma^1_0 = \#\Gamma^1_\infty - 3 + 3g$ bounded edges. A 3-valent graph of genus $g$ has $\#\Gamma^0_0 = \#\Gamma^1_\infty - 2 + 2g$ inner vertices. We need these relations for dimension counts later on.

*Remark 5.2* Consider an abstract tropical curve $\Gamma$ with $g$-labels. If $\Gamma$ is 3-valent and of genus $g$, then it has $s - 2 + 2g$ inner vertices, thus each inner vertex must be labeled with exactly one number. If $\Gamma$ has higher valent vertices or is of lower genus, then it is possible that $V_i = V_j$ for $i \neq j$ in this notation. We can think of such a curve as the result of shrinking several edges and $g - g'$ cycles in a 3-valent curve of genus $g$.

In the following, we define the tropical analogue of stable maps to $\mathbb{P}^1$. Note that tropical $\mathbb{P}^1$ can be thought of as $\mathbb{R} \cup \{\pm\infty\}$ (see e.g. [11]).

**Definition 5.3** Let $g \in \mathbb{N}$ and $\Delta$ be a multiset with entries in $\mathbb{Z} \setminus \{0\}$ satisfying $\sum_{z \in \Delta} z = 0$.

A (*parametrized*) *tropical curve of degree* $\Delta$ *and genus* $g$ *in* $\mathbb{P}^1$ is a tuple $(\Gamma, h)$ where $\Gamma$ is an abstract tropical curve with $g$-labels and $\#\Delta$ ends and $h : \Gamma \to \mathbb{R} \cup \{\pm\infty\}$ is a continuous map satisfying:

(a) The image of the graph without the 1-valent vertices has to be inside $\mathbb{R}$, $h(\Gamma \setminus \Gamma^0_\infty) \subset \mathbb{R}$.
(b) $h$ maps each edge $e$ of length $l(e)$ affinely to a line segment of $\mathbb{R} \cup \{\pm\infty\}$ of length $\omega(e) \cdot l(e)$, where $\omega(e)$ is a natural number that we call the *weight* of $e$.

For a flag $F = (V, e)$ we say that $F$ is of *direction* $v(F) := \omega(e)$ if $h(V) < h(p)$ for a point $V \neq p \in e$ and $F$ is of direction $v(F) := -\omega(e)$ otherwise. For ends $e$ we also say that the direction of $e$ is $v(e) := v(V, e)$, where $V$ is the inner end vertex of $e$.

(c) The multiset of directions of all ends equals $\Delta$.
(d) For every vertex $V \in \Gamma^0_0$ we have the *balancing condition*

$$\sum_{e | V \in \partial e} v(V, e) = 0.$$

Two parametrized tropical curves $(\Gamma, h)$ and $(\tilde{\Gamma}, \tilde{h})$ in $\mathbb{R}^r$ are called isomorphic (and will from now on be identified) if there is an isomorphism $\varphi : \Gamma \to \tilde{\Gamma}$ of the underlying abstract curves such that $\tilde{h} \circ \varphi = h$.

For the special choice $\Delta = \{-1, \dots, -1, 1 \dots, 1\}$ (each $d$ times) we also say that these curves have degree $d$.

*Remark 5.4* Note that 5.3 implies that a 1-valent vertex $V$ which is adjacent to an end $e$ satisfying $\omega(e) \neq 0$ has to be mapped to $\pm\infty$.

The *combinatorial type* of a tropical curve of degree $\Delta$ and genus $g$ in $\mathbb{P}^1$ is given by the data of the combinatorial type of the underlying abstract tropical curve $\Gamma$ together with the directions of all its flags (flags of bounded edges as well as ends).

The space of all tropical curves of degree $\Delta$ and genus $g$ in $\mathbb{P}^1$ is denoted $M_{g,\mathrm{trop}}(\mathbb{P}^1, \Delta)$.

*Remark 5.5* For a tropical curve $C = (\Gamma, h)$ and a point $p \in \mathbb{R}$ such that $h^{-1}(p)$ does not contain a vertex of $\Gamma$, the number of preimages (counted with the weight of the edge they are on) is constant because of the balancing condition (and equal to $d$ in case the degree is $d$). One can think of elements in $M_{g,\mathrm{trop}}(\mathbb{P}^1, \Delta)$ as limits of degree $d$ maps where we have "sent" some of the vertices to $\pm\infty$. I.e., we interpret the weight partitions over $\pm\infty$ (that are given by the positive/negative entries in $\Delta$) as special ramification profiles over the two points.

If we fix a combinatorial type $\alpha$ then we denote by $M^\alpha_{g,\mathrm{trop}}(\mathbb{P}^1, \Delta)$ the subset of curves in $M_{g,\mathrm{trop}}(\mathbb{P}^1, \Delta)$ with type $\alpha$. We call it a *cell* of $M_{g,\mathrm{trop}}(\mathbb{P}^1, \Delta)$.

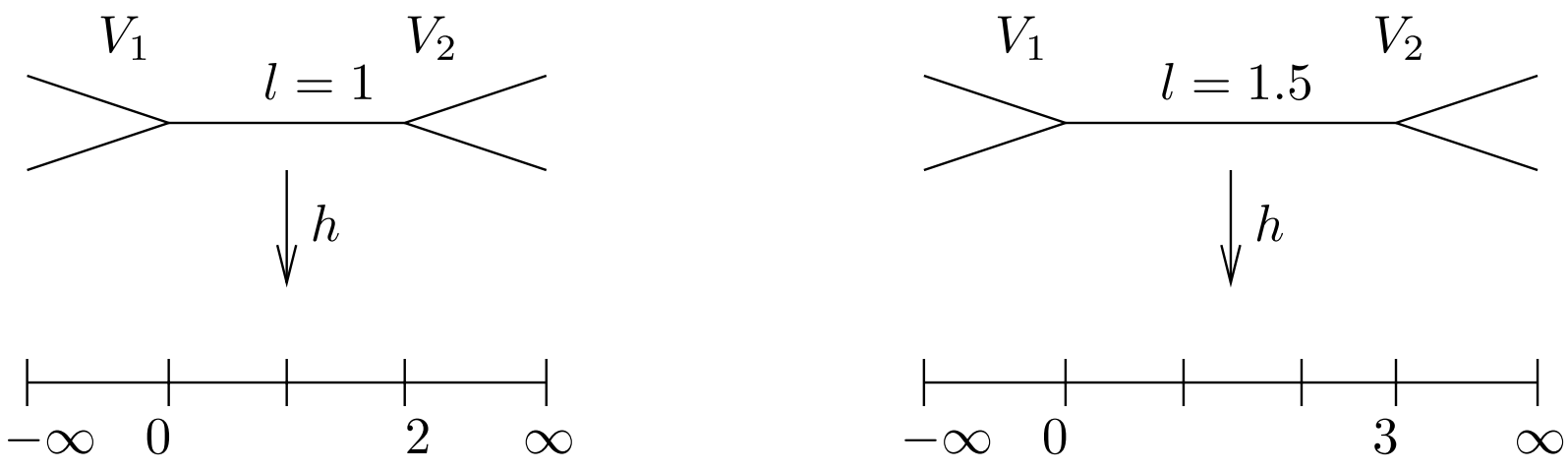

**Fig. 4** Two rational tropical curves in $\mathbb{P}^1$

*Example 5.6* Figure 4 shows two rational tropical curves of degree 2 in $\mathbb{P}^1$ of the same combinatorial type. If we denote the only bounded edge of the graph by $e$, then in both cases the direction of the flag $(V_1, e)$ is 2 and the direction of $(V_2, e)$ is $-2$. The multiset of directions of the four ends is $\{1, 1, -1, -1\}$.

**Lemma 5.7** *For a combinatorial type $\alpha$ of curves in $M_{g,\text{trop}}(\mathbb{P}^1, \Delta)$, the space $M^{\alpha}_{g,\text{trop}}(\mathbb{P}^1, \Delta)$ is an unbounded open convex polyhedron in a real vector space of dimension $1 + \#\Gamma^1_0$. It has one coordinate for the position of a root vertex and coordinates for the lengths of all bounded edges. $M^{\alpha}_{g,\text{trop}}(\mathbb{P}^1, \Delta)$ is cut out by the inequalities that all lengths have to be positive and by the equations for the loops.*

*The expected dimension of $M^{\alpha}_{g,\text{trop}}(\mathbb{P}^1, \Delta)$ is $1 + \#\Gamma^1_0 - g_\alpha = \#\Delta - 2 + 2g_\alpha - \sum_{V\in\Gamma^0_0}(\text{val}(V) - 3)$, where $g_\alpha \le g$ denotes the genus of $\Gamma$.*

The proof is a straightforward adaption of the proof of Lemma 3.1 of [8]. We say that a combinatorial type $\alpha$ is *regular* if $M^{\alpha}_{g,\text{trop}}(\mathbb{P}^1, \Delta)$ is of expected dimension.

Only types with a bounded flag of direction 0 can be non-regular. If we have a genus $g$ graph with no bounded flags of direction 0, we can pick $g$ edges such that $\Gamma$ without the $g$ edges is a tree. Then we can put back in one edge after the other. Each time, we close a loop and get a condition on the lengths of the edges in this loop, in particular, a condition on the edge we just put in. Since we put in the edges one after the other, we get a matrix with a triangular shape for the $g$ edges, thus independent conditions.

We want to make $M_{g,\text{trop}}(\mathbb{P}^1, \Delta)$ a *weighted polyhedral complex* of dimension $\#\Delta - 2 + 2g$ (which is the maximal expected dimension) in the sense of Definition 3.4 of [8].

In order to define weights for the maximal cells, we need the following notions: Let $f : \mathbb{Z}^n \to \mathbb{Z}^m$ be a linear map. We call the index of $f$, $I_f$, the index of the sublattice $f(\mathbb{Z}^n)$ inside $\mathbb{Z}^m$.

**Definition 5.8** Let $\alpha$ be a regular type of top dimension $\#\Delta - 2 + 2g$ in $M_{g,\text{trop}}(\mathbb{P}^1, \Delta)$. Pick $g$ independent cycles of the underlying 3-valent graph $\Gamma$, i.e. generators of $H_1(\Gamma, \mathbb{Z})$. Each such generator is given as a chain of flags around the loop. Define a $g$ times $1 + \#\Gamma^1_0 = \#\Delta - 2 + 3g$ matrix $A_\alpha$ with a column for the position of $h(V_1)$ and a column for each length coordinate, and with a row for each cycle containing

**Fig. 5** A combinatorial type

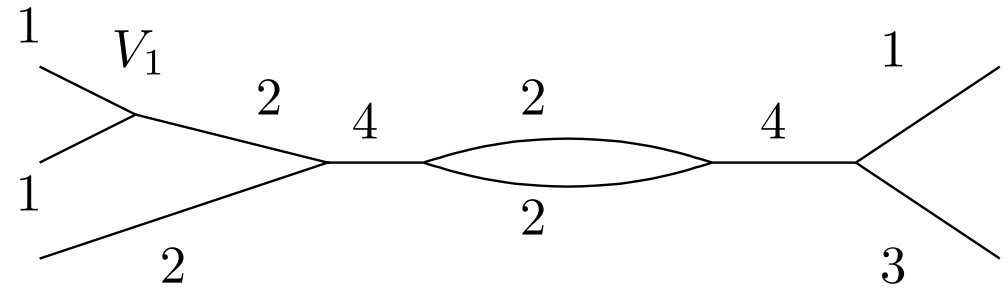

the equation of the loop (depending on the lengths of the bounded edges in the loop):

$$\sum_{(V,e)} v(V,e)\cdot l(e),$$

where the sum now goes over the chosen chain of flags around the loop. Denote by $I_\alpha$ the index of the map $A_\alpha : \mathbb{Z}^{\#\Delta-2+3g} \to \mathbb{Z}^g$. Note that $I_\alpha$ does not depend on the chosen generators of $H_1(\Gamma,\mathbb{Z})$: if we choose another set of generators, these new generators are given as linear combinations with coefficients in $\mathbb{Z}$ of the old generators, so the rowspace of the matrix is not changed. Note also that $M^\alpha_{g,\mathrm{trop}}(\mathbb{P}^1,\Delta)$ equals the intersection of $\mathbb{R}\times(\mathbb{R}_{>0})^{\#\Gamma^1_0}$ with the kernel of this map.

*Example 5.9* For a curve of type $\alpha$ as in Fig. 5, we have 5 bounded edges and one loop. $A_\alpha$ is a $1\times 6$ matrix. The loop is formed by two edges, say $e_1$ and $e_2$. Going around the loop, we set up the equation $2\cdot l(e_1)+(-2)\cdot l(e_2)$, so the matrix reads $(0,2,-2,0,0,0)$. The entries in the column of the root vertex are always 0.

We have to throw away cells of $M_{g,\mathrm{trop}}(\mathbb{P}^1,\Delta)$ of too big dimension. When we introduce the branch map later and define its degree, cells of too big dimension would not contribute to the count anyway. However, we would like the dimension of the space as expected.

**Definition 5.10** Let $\mathcal{M}_{g,\mathrm{trop}}(\mathbb{P}^1,\Delta)$ be the subset of $M_{g,\mathrm{trop}}(\mathbb{P}^1,\Delta)$ containing all combinatorial types $\alpha$ such that if $M^\alpha_{g,\mathrm{trop}}(\mathbb{P}^1,\Delta)$ is of dimension $\#\Delta-2+2g$ or bigger then $\alpha$ is regular and if $M^\alpha_{g,\mathrm{trop}}(\mathbb{P}^1,\Delta)$ is of dimension less than $\#\Delta-2+2g$ then it is contained in a cell $M^{\alpha'}_{g,\mathrm{trop}}(\mathbb{P}^1,\Delta)$ of highest dimension. In particular, the highest dimension of a cell in $\mathcal{M}_{g,\mathrm{trop}}(\mathbb{P}^1,\Delta)$ is $\#\Delta-2+2g$.

Let $\alpha$ be a type corresponding to a cell of highest dimension. Recall that the data of $\alpha$ consists of the graph $\Gamma$ (without lengths) and the information about all directions. We define its *weight* $w(\alpha)$ as the product of three types of factors:

- $\frac{1}{2}$ for every vertex $V$ such that $\Gamma\setminus V$ has two connected components of the same combinatorial type.
- The index $I_\alpha$.
- $\frac{1}{2}$ for every cycle which consists of two edges which have the same weight, i.e. for every wiener as in Lemma 4.2.

*Remark 5.11* These choices of weights are not new to tropical geometry (see e.g. [8], Remark 3.6 or [4]). We can interpret the factors of $\frac{1}{2}$ as taking care of extra automorphisms.

*Example 5.12* Figure 5 shows a type of curves of genus 1 and degree $\{-1,-1,-2,1,3\}$. The corresponding cell gets the weight $\frac{1}{2}\cdot\frac{1}{2}\cdot 2$. The first factor of $\frac{1}{2}$ appears because two of the connected components of $\Gamma \setminus \{V_1\}$ are identical. The second factor appears because of the wiener. The third factor is the index of the map $(0,2,-2,0,0,0)$ which is the equation for the loop. The numbers written next to the edges are the weights of the edges (there are no lengths, because the picture only shows a combinatorial type).

Finally we have to glue the cells $M^{\alpha}_{g,\mathrm{trop}}(\mathbb{P}^1,\Delta)$ to make $\mathcal{M}_{g,\mathrm{trop}}(\mathbb{P}^1,\Delta)$ a weighted polyhedral complex. This can be done analogously to Proposition 3.2 of [8].

**Lemma 5.13** *The tropical moduli space $\mathcal{M}_{g,\mathrm{trop}}(\mathbb{P}^1,\Delta)$ as defined in Definition 5.10 is a weighted polyhedral complex of pure dimension $\#\Delta-2+2g$.*

5.2 The tropical branch map

**Definition 5.14** We define the tropical branch map *br* as

$$br:\mathcal{M}_{g,\mathrm{trop}}\big(\mathbb{P}^1,\Delta\big)\to\mathbb{R}^{\#\Delta-2+2g}:(\Gamma,h)\mapsto\big(h(V_1),\dots,h(V_{\#\Delta-2+2g})\big),$$

where $V_i$ is the vertex or point with label $i$.

*Example 5.15* For the left picture of Example 5.6, we have $br(\Gamma,h)=(h(V_1),h(V_2))=(0,2)$.

Then *br* is a *morphism of weighted polyhedral complexes* of the same dimension in the sense of Definition 4.1 of [8]. To see this, we have to see that it is a linear map on each cell $M^{\alpha}_{g,\mathrm{trop}}(\mathbb{P}^1,\Delta)$. This is true because the position $h(V_i)$ differs from $h(V_1)$ by a sum $v(V_1,e_1)l(e_1)+\cdots+v(V_r,e_r)l(e_r)$ where $(V_1,e_1),\dots,(V_r,e_r)$ denotes a chain of flags that we have to pass to go from $V_1$ to $V_i$ in $\Gamma$. Note that the map does not depend on the chain of flags we choose: going another way around a cycle does not change anything since the length coordinates satisfy the conditions given by the cycles.

**Definition 5.16** For a type $\alpha$ of maximal dimension in $\mathcal{M}_{g,\mathrm{trop}}(\mathbb{P}^1,\Delta)$ we choose the following data:

- for each vertex $V_i\in\Gamma^0_0$ a chain of flags leading from $V_1$ to $V_i$ and
- a set of generators of $H_1(\Gamma,\mathbb{Z})$, where each such generator is given as a chain of flags around the loop.

Depending on these choices, we define a linear map $f_\alpha$ by defining a square matrix of size $\#\Delta-2+3g$ with

- for each vertex $V_i\in\Gamma^0_0$ a row with the linear equation describing the position of $h(V_i)$ (depending on the position of $h(V_1)$ and the lengths of the bounded edges in the chosen chain of flags from $V_1$ to $V_i$):

$$h(V_1)+\sum_{(V,e)}v(V,e)\cdot l(e),$$

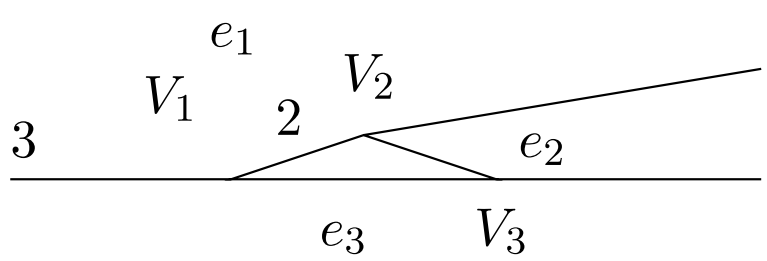

**Fig. 6** The curve of Example 5.18

where the summation goes over all flags $(V, e)$ in the chosen chain from $V_1$ to $V_i$; and

- $g$ rows as in the matrix $A_\alpha$ defined in Definition 5.8.

*Remark 5.17* The map $f_\alpha$ of Definition 5.16 depends on the choices we made, while the absolute value of the determinant of $f_\alpha$ does not. First, a coordinate change for $M^\alpha_{g,0,\mathrm{trop}}(\mathbb{P}^1, \Delta)$ has determinant $\pm 1$ and therefore leaves the absolute value of $\det f_\alpha$ unchanged. The same holds for a coordinate change of the target space $\mathbb{R}^{\#\Delta-2+3g}$, that is, a different order of the vertices and loops.

If there are two chains of flags from $V_1$ to $V_i$, then their difference is a loop. Assume first that this loop is one of our chosen generators of $H_1(\Gamma, \mathbb{Z})$. Then choosing one or the other chain of flags from above just corresponds to adding (respectively subtracting) the equation of the loop from the row of $V_i$. We have seen already in Definition 5.8 that other generators $H_1(\Gamma, \mathbb{Z})$ are given as linear combinations with coefficients in $\mathbb{Z}$ of the old generators.

By abuse of notation, we still speak of the map $f_\alpha$, even though its definition depends on the choices we made, and keep in mind that $|\det(f_\alpha)|$ is uniquely determined, no matter what choices we made.

*Example 5.18* For the curve in Fig. 6, we have 3 vertices, $V_1$, $V_2$ and $V_3$.

As usual, $V_1$ is the root vertex. We choose to go to $V_2$ via $e_1$ and to $V_3$ via $e_1$ and $e_2$. So $h(V_2) = h(V_1) + 2 \cdot l(e_1)$ and $h(V_3) = h(V_1) + 2 \cdot l(e_1) + 1 \cdot l(e_2)$. There is one loop which consists of the three edges $e_1$, $e_2$ and $e_3$. So the equation for the loop is $2 \cdot l(e_1) + 1 \cdot l(e_2) - 1 \cdot l(e_3)$. That amounts to the following matrix for $f_\alpha$:

$$\begin{pmatrix} 1 & 0 & 0 & 0 \\ 1 & 2 & 0 & 0 \\ 1 & 2 & 1 & 0 \\ 0 & 2 & 1 & -1 \end{pmatrix}.$$

The absolute value of its determinant is 2.

*Remark 5.19* A straightforward lattice index computation shows that for a combinatorial type $\alpha$ of maximal dimension in $\mathcal{M}_{g,\mathrm{trop}}(\mathbb{P}^1, \Delta)$, the index $I_\alpha$ times the absolute value of the determinant of the linear map $br$ restricted to the cell $M^\alpha_{g,\mathrm{trop}}(\mathbb{P}^1, \Delta)$ is equal to $|\det(f_\alpha)|$. A similar lattice index computation can e.g. be found in Remark 4.8 in [8], or see [12], Lemma 1.6. The index of a square integer matrix is just the absolute value of its determinant, and the index of a product of two maps $f \times g$ is equal to the index of $f|_{\ker g}$ times the index of $g$.

**Fig. 7** An equivalence class of combinatorial types

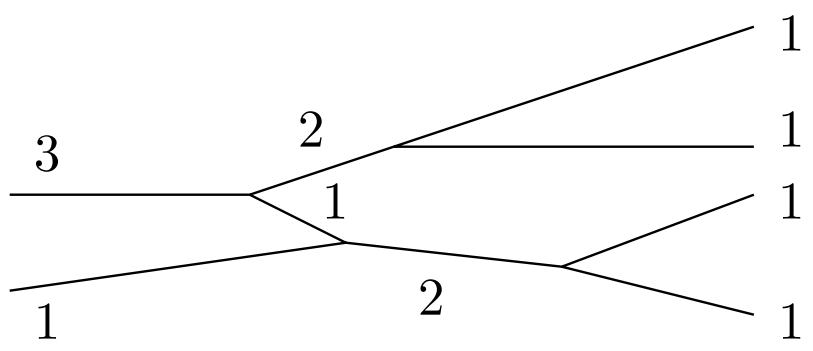

**Definition 5.20** Two combinatorial types $\alpha$ and $\alpha'$ are called *equivalent* if they differ only by the labeling of the vertices.

*Example 5.21* We can visualize an equivalence class of combinatorial types by a graph without labeling for the vertices, but together with the information about the weights of the edges. In this sense, Fig. 7 shows an equivalence class of combinatorial types of rational curves of degree $\{-3, -1, 1, 1, 1, 1\}$.

*Remark 5.22* The set of points $p \in \mathbb{R}^{\#\Delta-2+2g}$ satisfying $p_i \neq p_j$ is an open dense subset of $\mathbb{R}^{\#\Delta-2+2g}$. It is contained in the set of points in *br*-general position, because for any preimage under *br*, no vertex has more than one number as label and so each vertex has to be 3-valent. Then the expected dimension of the combinatorial type has the highest dimension. Thus the type is regular and of expected dimension, since we removed types of too high dimension.

**Definition 5.23** Let $[\alpha]$ be an equivalence class of types of highest dimension in $\mathcal{M}_{g,\mathrm{trop}}(\mathbb{P}^1, \Delta)$. We define a partial ordering on the vertices in $\Gamma_0^0$ in the following way: $V < V'$ if and only if $V'$ can be reached from $V$ by a chain of flags with positive direction. We denote by $n([\alpha])$ the number of ways to extend this ordering to a well-ordering.

**Lemma 5.24** *Fix a class $[\alpha]$ of combinatorial types of highest dimension in $\mathcal{M}_{g,\mathrm{trop}}(\mathbb{P}^1, \Delta)$. Fix $p \in \mathbb{R}^{\#\Delta-2+2g}$ satisfying $p_i \neq p_j$. Then there are $n([\alpha])$ combinatorial types $\alpha'$ in the class of $\alpha$ such that there is a preimage of $p$ under br in $M'^{\alpha'}_{g,0,\mathrm{trop}}(\mathbb{P}^1, \Delta)$.*

*Proof* Without loss of generality, we can assume $p_1 < \cdots < p_{\#\Delta-2+2g}$. If $V < V'$ in the partial ordering on the vertices of $\Gamma$ then the images $h(V)$ and $h(V')$ have to satisfy $h(V) < h(V')$. If we choose one of the $n([\alpha])$ well-orderings extending $<$, then there is only one labeling of the vertices which can satisfy $h(V_i) = p_i$. Let $\alpha'$ be the combinatorial type of class $[\alpha]$ with this labeling. We have to show that there is a tropical curve of this type in the preimage of $p$ under $br$. $\Gamma$ has $\#\Delta - 3 + 3g$ bounded edges. Each edge has two end vertices. Therefore the length of the edge is given by the distance of the images $h(V_i) = p_i$ of the two end vertices. We only need to show that those lengths satisfy the conditions given by the loops. But this is obviously true, since the images of the vertices contained in a loop close up a loop in $\mathbb{R}$. □

*Example 5.25* For the equivalence class of types in Example 5.21, there are 3 different choices for a well-ordering on the vertices extending the ordering of Lemma 5.24 as shown in Fig. 8.

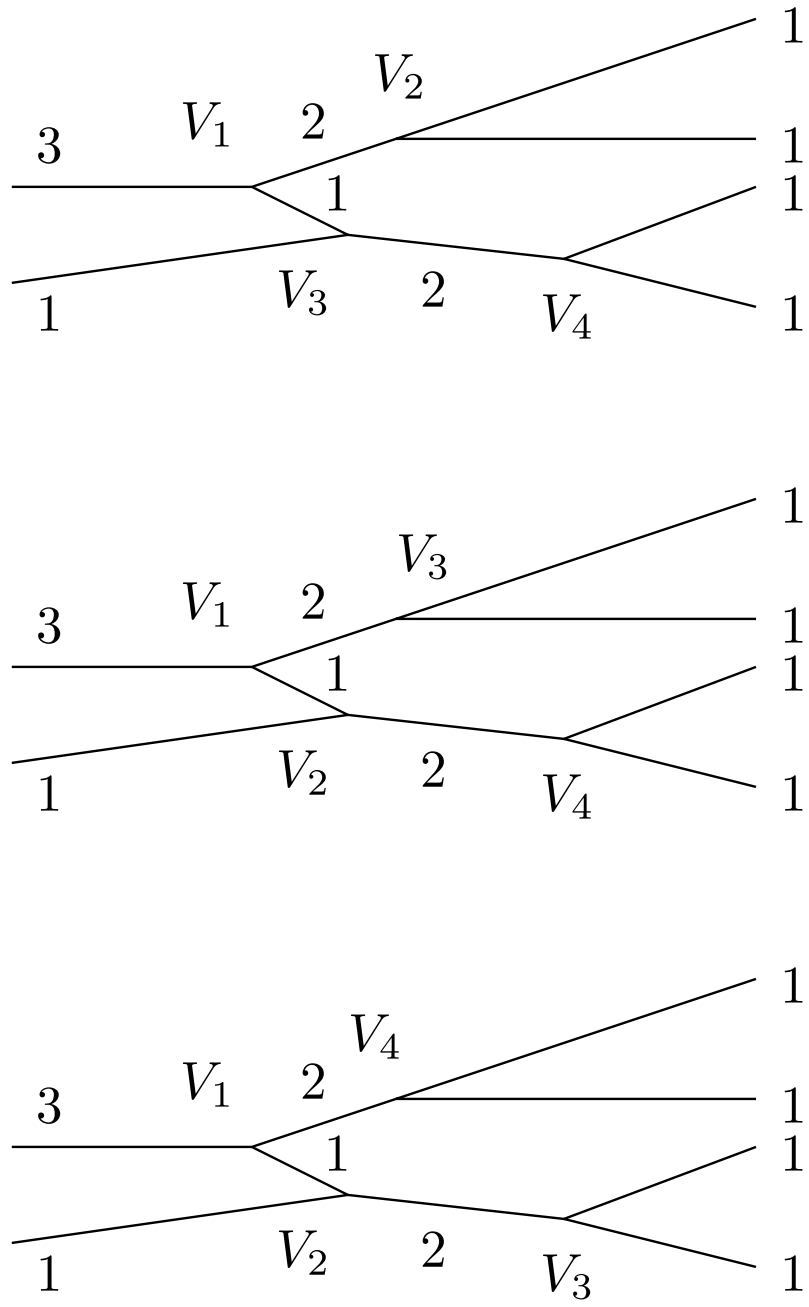

**Fig. 8** Ordering the vertices

**Lemma 5.26** *Fix an equivalence class* $[\alpha]$ *of types in* $\mathcal{M}_{g,\mathrm{trop}}(\mathbb{P}^1, \Delta)$ *of the highest dimension and* $p \in \mathbb{R}^{\#\Delta-2+2g}$ *satisfying* $p_i \neq p_j$. *Then the contribution* $k_{[\alpha]}$ *of curves of types in this class to* $\deg_{br}(p)$ *is given by*

$$k_{[\alpha]} = n\big([\alpha]\big) \cdot \left(\frac{1}{2}\right)^{r+s} \cdot \prod_e \omega(e), \tag{3}$$

*where* $r$ *denotes the number of vertices* $V$ *such that* $\Gamma \setminus V$ *has two identical parts*, $s$ *denotes the number of wieners* (*as in Lemma* 4.2) *and the product goes over all bounded edges* $e$ *in* $\Gamma_0^1$.

*Proof* Let $\alpha$ and $\alpha'$ be types in the class $[\alpha]$. Obviously $\mathrm{mult}_{br}(C) = \mathrm{mult}_{br}(C')$ for two curves $C \in M^{\alpha}_{g,\mathrm{trop}}(\mathbb{P}^1, \Delta)$ and $C' \in M'^{\alpha'}_{g,0,\mathrm{trop}}(\mathbb{P}^1, \Delta)$. Therefore the contribution of types of class $[\alpha]$ is equal to $k_{[\alpha]} = n([\alpha]) \cdot \mathrm{mult}_{br}(C)$ (where $C \in M^{\alpha}_{g,\mathrm{trop}}(\mathbb{P}^1, \Delta)$). The latter multiplicity is by definition equal to the weight of $\alpha$, $w(\alpha)$, times the absolute value of the determinant of the linear map $br$ restricted to $M^{\alpha}_{g,\mathrm{trop}}(\mathbb{P}^1, \Delta)$. By definition, $w(\alpha) = I_\alpha \cdot (\frac{1}{2})^{r+s}$ and by Remark 5.19, the absolute value of the determinant of $br$ times $I_\alpha$ equals $|\det(f_\alpha)|$. Thus $k_{[\alpha]} = n([\alpha]) \cdot (\frac{1}{2})^{r+s} \cdot |\det(f_\alpha)|$ and it remains to show that $|\det(f_\alpha)| = \prod_e \omega(e)$.

To see this, we want to show that we can choose an order of the coordinates such that the matrix of $f_\alpha$ is lower triangular and the weights of all bounded edges appear on the diagonal. Start by removing $g$ bounded edges from $\Gamma$ breaking the cycles. The new graph that we call $\Gamma'$ is rational and has the same set of vertices. It has 1-, 2- and 3-valent vertices in $\Gamma_0^0$. We want to show that there is an order of the edges and of the vertices in $\Gamma_0^0$ such that we need only edges of order less than $i$ to go from vertex 1 to vertex $i$, and such that edge $i-1$ is adjacent to vertex $i$ and needed in the path from $V_1$ to $V_i$. We show this by induction. The induction beginning—where $\Gamma'$ has only one vertex in $\Gamma_0^0$—is obvious. Now we can assume $\Gamma'$ has at least two vertices in $\Gamma_0^0$. The subgraph of $\Gamma'$ of bounded edges is a tree, so it has to have at least two leaves. This means that there are at least two vertices which are adjacent to only one bounded edge. Call one of those vertices $V_k$ and its adjacent bounded edge $e_{k-1}$. Remove $V_k$ and all adjacent ends and make $e_{k-1}$ an end. We call this new graph $\Gamma''$. By induction we can assume that we can order the remaining edges and vertices in $\Gamma_0^0$ in the way we require. Now we add $V_k$ and $e_{k-1}$ back in. To go from $V_1$ to $V_k$, we need edges in $\Gamma''$ and $e_{k-1}$. Now add the $g$ edges back in, one after the other. At each step, we close a loop. We write down the equation for this loop as the next row of our matrix. Thus we end up with a lower triangular matrix such that the absolute value of its determinant is equal to the product of all weights of bounded edges. □

**Lemma 5.27** *The degree of br is constant, i.e.* $\deg(br) := \deg_{br}(p)$ *does not depend on the choice of* $p$, *as long as we pick* $p$ *such that* $p_i \neq p_j$.

*Proof* By the above, $\deg_{br}(p) = \sum k_{[\alpha]}$ where the sum goes over all equivalence classes of combinatorial types of highest dimension in $\mathcal{M}_{g,\mathrm{trop}}(\mathbb{P}^1, \Delta)$. □

### 5.3 The main theorem

We combine the results from Sects. 4 and 5.2 to prove our main theorem:

**Theorem 5.28** *The degree* $\deg(br)$ *of the tropical branch map* $br : \mathcal{M}_{g,\mathrm{trop}}(\mathbb{P}^1, \Delta) \to \mathbb{R}^{\#\Delta-2+2g}$ *as defined in Definition 5.14 is equal to the Hurwitz number* $H_d^g(\eta, \nu)$ *as defined in Definition 2.1, where* $\eta$ *is the partition of* $d$ *given by the negative entries in* $\Delta$ *and* $\nu$ *is the partition of* $d$ *given by the positive entries in* $\Delta$.

*Proof* From Lemmas 5.27 and 4.2 we know that the degrees of the tropical branch map and the Hurwitz number $H_g^d(\eta, \nu)$ are computed via a weighted sum of graphs. We now show that, for each combinatorial type of graphs, we have the same coefficient in both cases. Fix a combinatorial type $\alpha$.

In Lemma 5.26 the contribution from graphs of type $\alpha$ is given by (3). After simplifying via (1) in Corollary 4.4, the contribution on the Hurwitz count side is

$$n'([\alpha]) \cdot \left(\frac{1}{2}\right)^{\#\text{b.forks}+\#\text{wieners}} \cdot \prod_e \omega(e), \tag{4}$$

where $n'([\alpha])$ is the number of times a graph of type $\alpha$ appears in the construction, and the product is over all bounded edges. To prove the theorem it suffices to show:

$$n'([\alpha]) \cdot \left(\frac{1}{2}\right)^{\#\text{b.forks}} = n([\alpha]) \cdot \left(\frac{1}{2}\right)^{r},$$

where $r$ is the number of vertices $V$ such that $\Gamma \setminus V$ has two identical connected components.

Graphs in the Hurwitz count induce a labeling of the vertices by assigning label $i$ to the vertex mapping to $i$. However we are counting isomorphism classes of labeled graphs: a graph has two indistinguishable labelings precisely when removing a vertex (non-adjacent to any end) there are two identical connected components, i.e.

$$n'\big([\alpha]\big) = n\big([\alpha]\big) \cdot \left(\frac{1}{2}\right)^{\text{r}-\#\text{b.forks}},$$

which is what we need to show. □

*Remark 5.29* A geometric proof of our main theorem should also follow from Theorem 1, [11]. We thank G. Mikhalkin for pointing this out to us.

## 6 Combinatorial properties of genus 0 double Hurwitz numbers

In this section we use Lemma 4.2 to recover in an elementary way results in [6] and [13] on the structure of double Hurwitz numbers in genus 0. That is, we do not view Hurwitz numbers one by one, but as a function on the entries of the two partitions $\mu$ and $\nu$. We point the attention of the reader to a technical detail: in this section we wish to adopt the definition of Hurwitz numbers in [6], which differs from the classical one used so far in that the preimages of 0 and $\infty$ are marked. The difference between the two definitions is a multiplicative factor of $|\mathrm{Aut}(\mu)||\mathrm{Aut}(\nu)|$.

**Definition 6.1** Let $k+l \geq 3$ and $\mu_1, \dots, \mu_k, \nu_1, \dots, \nu_l$ be the coordinates of $\mathbb{R}^{k+l}$. Let $\mathcal{H}$ be the hyperplane $\sum \mu_i = \sum \nu_j$. We think of $H^0$ as a map

$$\begin{aligned} &H^0 : \mathcal{H} \cap \mathbb{N}^{k+l} \\ &\quad \to \mathbb{Q} : \quad (\mu_1, \dots, \mu_k, \nu_1, \dots, \nu_l) \mapsto H^0_{\mu_1+\dots+\mu_k}\big((\mu_1, \dots, \mu_k), (\nu_1, \dots, \nu_l)\big). \end{aligned}$$

### 6.1 Piecewise polynomiality

**Theorem 6.2** [6, 13] *The map $H^0$ is piece-wise polynomial. More precisely, $\mathcal{H}$ is subdivided into a finite number of chambers, and inside each chamber the map $H^0$ is a homogeneous polynomial in the $\mu_i$ and $\nu_j$ of degree $k+l-3$. Walls defining the chambers are given by the equations*:

$$\sum_{i \in I} \mu_i - \sum_{j \in J} \nu_j = 0, \tag{5}$$

*for $I$, $J$ any proper subsets of the indices sets.*

We introduce some notions that are necessary to the proof of Theorem 6.2.

**Definition 6.3** Let $\mathcal{T}(k,l)$ be the set of all connected 3-valent directed trees with $k$ in-ends labeled with $1,\dots,k$ and $l$ out-ends labeled with $1,\dots,l$ and with no sources or sinks, together with a total order of the vertices compatible with the edge directions.

If we assign weights to the ends of any graph in $\mathcal{T}(k,l)$ (in such a way that the sum of the weights of the in-ends equals the sum of the weights of the out-ends), we can then weight the internal edges in a unique way by imposing the balancing condition.

**Lemma 6.4** *Let $\Gamma \in \mathcal{T}(k,l)$. If we choose the weight $\mu_i$ for the in-end labeled $i$ and $\nu_j$ for the out-end labeled $j$, then the weights $\omega(e)$ of all inner edges $e$ are uniquely determined* (*but might be negative*). *The weight $\omega(e)$ equals*

$$\omega(e) = \sum_{i\in I}\mu_i - \sum_{j\in J}\nu_j, \tag{6}$$

*where $I \subset \{1,\dots,k\}$ and $J \subset \{1,\dots,l\}$ are the subsets of in- and out-ends belonging to the connected component of $\Gamma \setminus \{e\}$ from which $e$ points away.*

*Proof* The balancing condition implies that the weighted graphs can be interpreted as networks of flowing water where water is neither created nor destroyed: the sum of inflow and outflow must then be equal. When cutting a tree along an internal edge, the two resulting connected components also satisfy this condition, thus determining uniquely the weight of the cut edge to be given by formula (6). □

We note that the weights of all internal edges are linear homogeneous polynomials in the entries of the partitions.

**Definition 6.5** For a fixed pair of partitions $\mu,\nu$, we define $\mathcal{T}^+(\mu,\nu) \subset \mathcal{T}(k,l)$ to be the subset of graphs such that all internal edges have strictly positive weights when the ends are given weights corresponding to the partitions $\mu$ and $\nu$.

In formula (2), we are summing over all graphs $\Gamma \in \mathcal{T}^+(\mu,\nu)$: for every $\Gamma \in \mathcal{T}^+(\mu,\nu)$, we can build a projection sending $\Gamma$ to the interval $[0,s+1]$ that maps the source vertices of the in-ends to 0, the target vertices of the out-ends to $s+1$ and the vertex with label $i$ to $i$. The projection satisfies that the image of the source vertex of an edge is smaller than the image of the target vertex. Vice versa, by directing the edges in a graph projecting to $[0,s+1]$, we get an element in $\mathcal{T}^+(\mu,\nu)$. From Lemma 6.4 we see that the set $\mathcal{T}^+(\mu,\nu)$ is constant precisely in the chambers $C$ defined by the walls in (5); thus we also use $\mathcal{T}^+(C)$ to denote $\mathcal{T}^+(\mu,\nu)$ for any $\mu,\nu$ in $C$. We have now proved all ingredients needed in the proof of Theorem 6.2.

*Proof of Theorem 6.2* Proving this theorem using formula (2) is elementary: for any contributing graph the weights of the internal edges are linear homogeneous polynomials in the $\mu_i$'s and $\nu_j$'s (Lemma 6.4), and we are taking a product over $k+l-3$ internal edges. We next sum over a finite set of graphs, which remains constant in regions where the signs of all internal edges does not change (Definition 6.5). □

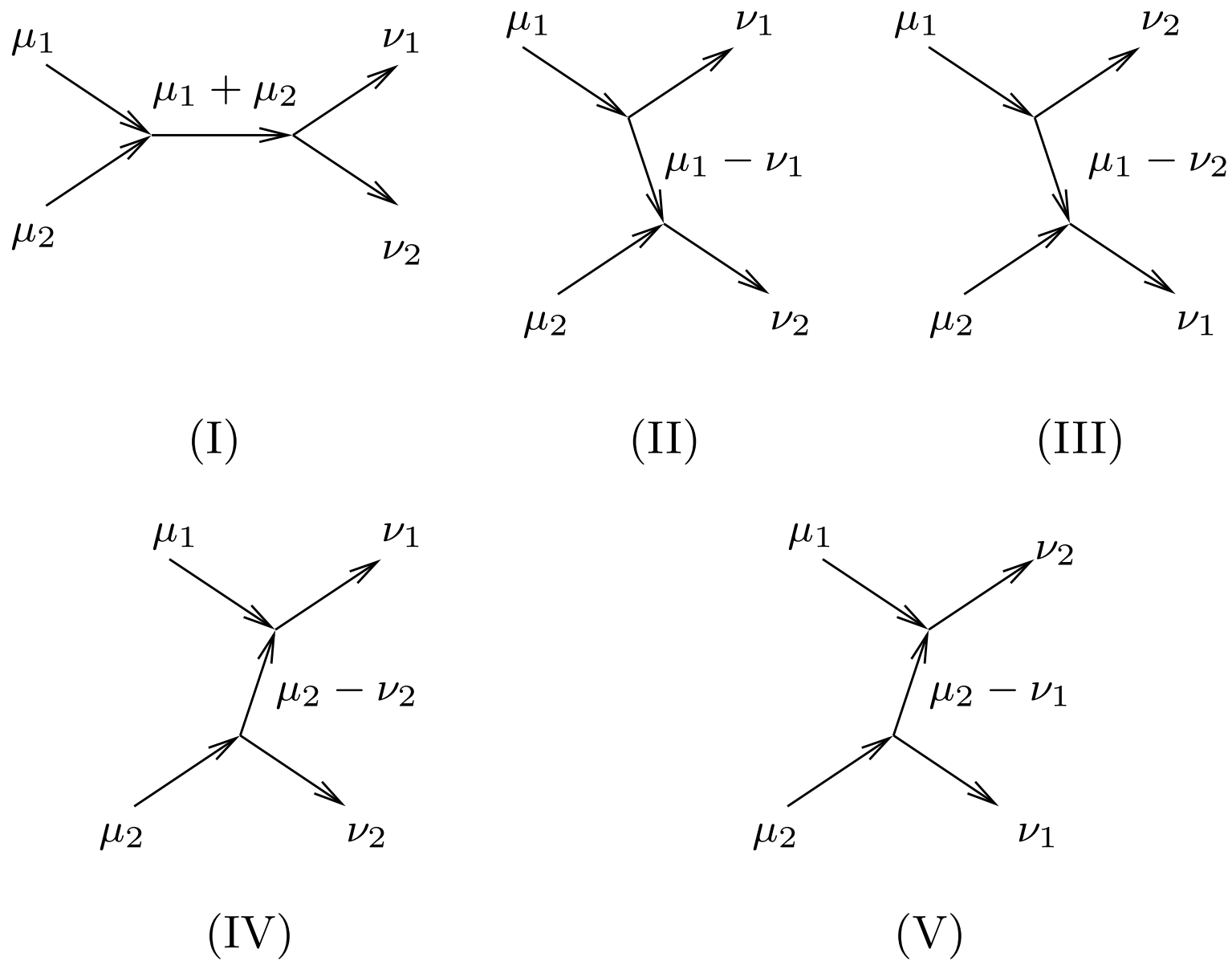

**Fig. 9** Weighted trees with two in-ends and two out-ends

*Example 6.6* Figure 9 shows graphs of the set $\mathcal{T}(2,2)$, after attaching the weights for the ends and concluding the weight of the inner edge according to Lemma 6.4.

The chambers are defined by the inequalities $\pm(\mu_1-\nu_1)>0$, and $\pm(\mu_1-\nu_2)>0$. Note that two such inequalities, e.g. $\mu_1>\nu_1$ and $\mu_1>\nu_2$ imply two other inequalities $\nu_1>\mu_2$ and $\nu_2>\mu_2$. This is true since the sum $\mu_1+\mu_2$ equals $\nu_1+\nu_2$. Figure 10 shows the four chambers and the two walls, and marks which of the above graphs belongs to $\mathcal{T}^+(C)$ for each chamber $C$. Also, it shows the polynomial which equals $H^0$ in each chamber.

## 6.2 Wall crossing formulas

In this section we investigate how the polynomials computing Hurwitz numbers vary from chamber to chamber. Example 6.6 suggests a crucial observation: if we are only concerned with the difference of the polynomials across a wall $\delta=0$, we need only consider the contributions from graphs that belong to $\mathcal{T}^+$ in only one of the two chambers in questions. Further we can characterize these graphs as those containing an edge with weight $\delta$ that switches direction across the wall. This allows us to easily recover the formulas of [13].

**Definition 6.7** Choose a subset $I$ of the in-ends, and a subset $J$ of the out-ends. This defines a wall $\delta=\sum_{i\in I}\mu_i-\sum_{j\in J}\nu_j=0$. We select two adjacent chambers $C_1$ and $C_2$: all inequalities defining these chambers are the same, except for the inequality corresponding to the wall. We define $C_1$ to correspond to $\sum_{i\in I}\mu_i-\sum_{j\in J}\nu_j>0$. Let $P_1$ denote the polynomial that equals $H^0$ in $C_1$, $P_2$ the polynomial that equals $H^0$ in $C_2$. Then the wall crossing for this wall and the two adjacent chambers $C_1$ and $C_2$ is defined to be the difference of the two polynomials:

$$\mathrm{WC}_\delta(\mu,\nu) := P_1(\mu,\nu) - P_2(\mu,\nu).$$

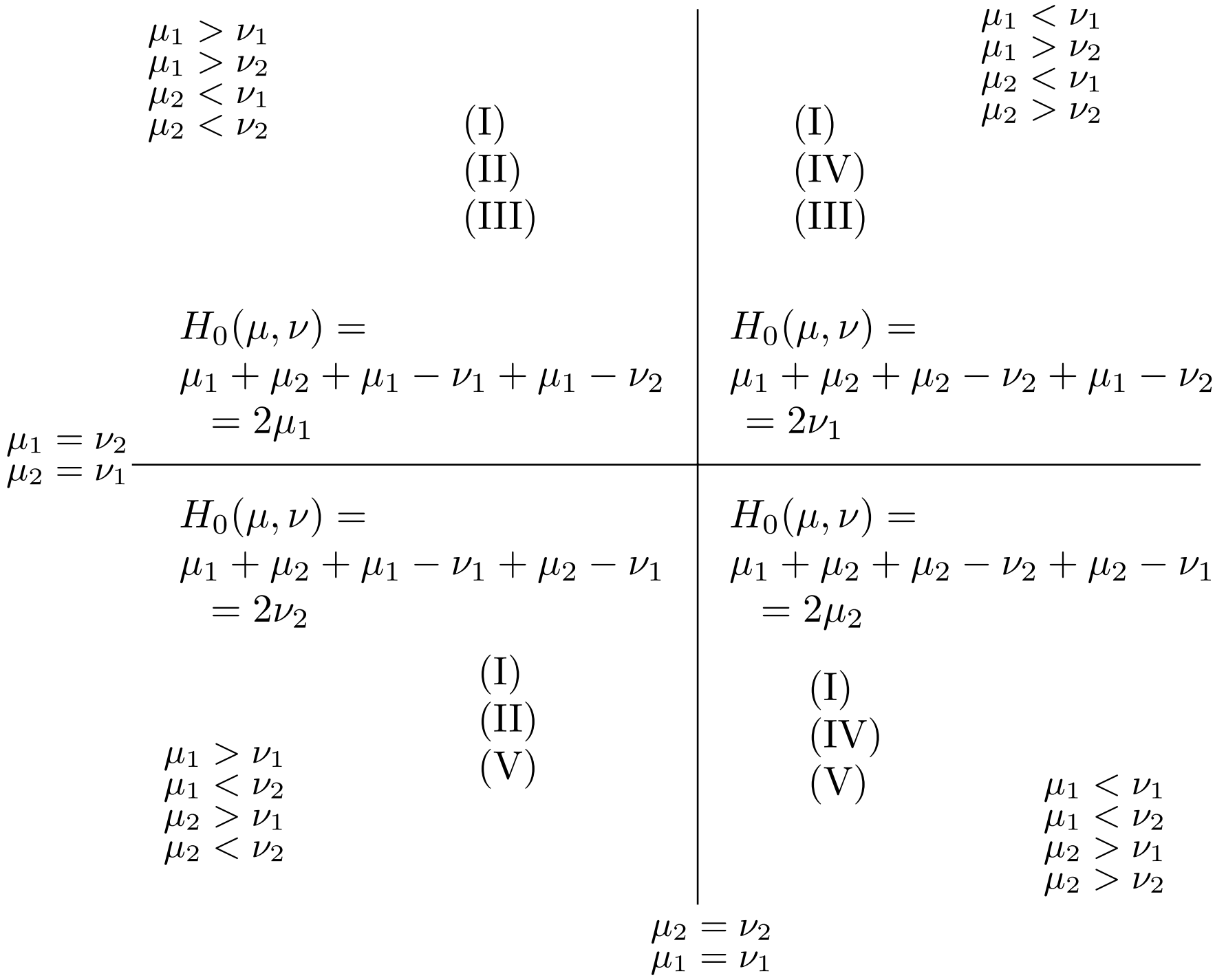

**Fig. 10** The polynomiality chambers for Example 6.6

*Example 6.8* Refer to Example 6.6 (Fig. 10) and consider the wall $\mu_1 = \nu_1, \mu_2 = \nu_2$ and the two adjacent regions on the top, i.e. $\mu_1 > \nu_1, \mu_1 > \nu_2, \mu_2 < \nu_1, \mu_2 < \nu_2$ and $\mu_1 < \nu_1, \mu_1 > \nu_2, \mu_2 < \nu_1, \mu_2 > \nu_2$. The graph (II) appears on the left, but is replaced on the right by (IV). This happens because this graph has an edge with weight going to 0 when approaching the wall. Graph (IV) equals graph (II) with the direction of the edge going to 0 reversed, and the weight multiplied with $-1$ accordingly.

The graphs (I) and (III) appear on both sides of the wall. Their contribution to the respective polynomials remains unchanged and they thus do not contribute to the wall crossing.

Thus, the wall crossing is given by the contribution from (II) minus the contribution from (IV), i.e. $\mu_1 - \nu_1 - (-\mu_1 + \nu_1) = 2\mu_1 - 2\nu_1$.

**Lemma 6.9** *For the wall* $\delta = 0$ *the wall crossing* $\mathrm{WC}_\delta(\mu,\nu)$ *can be expressed in terms of graphs in* $\mathcal{T}^\delta := (\mathcal{T}^+(C_1) \cup \mathcal{T}^+(C_2)) \setminus (\mathcal{T}^+(C_1) \cap \mathcal{T}^+(C_2))$; *all these graphs* *have an edge $\bar{e}$ with weight $\pm\delta$:*

$$\mathrm{WC}_\delta(\mu,\nu) = \sum_{\Gamma \in \mathcal{T}^\delta} \delta \cdot \prod_{e' \neq \bar{e}} \omega(e').$$

*Proof* We have

$$\mathrm{WC}_\delta = P_1 - P_2 = \sum_{\Gamma \in \mathcal{T}^+(C_1)} \varphi_\Gamma - \sum_{\Gamma \in \mathcal{T}^+(C_2)} \varphi_\Gamma,$$

where $\varphi_\Gamma$ denotes the product of the weights of all internal edges of $\Gamma$, as in formula (2) (note that the graphs have no automorphisms: since the genus is 0 there are no wieners, and all ends are labeled with distinct variables). Both chambers are given by inequalities of the form $\sum_{i\in I'} \mu_i - \sum_{j \in J'} \nu_j > 0$. The only inequality that differs for the two regions is the equality of the wall: $\delta = 0$. Let us say in $C_1$ we have $\delta > 0$ and in $C_2$ we have $\delta < 0$.

If $\Gamma$ is a graph in $\mathcal{T}^+(C_1)$ that does not have an edge of weight $\delta$, then all the weights will be positive in both chambers since all other inequalities persist. Thus $\Gamma$ also belongs to $\mathcal{T}^+(C_2)$ and hence its contribution $\varphi_\Gamma$ cancels in the difference $P_1 - P_2$.

We then have

$$\mathrm{WC}_\delta = P_1 - P_2 = \sum_{\Gamma \in (\mathcal{T}^+(C_1) \setminus \mathcal{T}^+(C_2))} \varphi_\Gamma - \sum_{\Gamma \in (\mathcal{T}^+(C_2) \setminus \mathcal{T}^+(C_1))} \varphi_\Gamma.$$

Notice that for each graph in $\mathcal{T}^\delta$ there is precisely one edge of weight $\pm\delta$. This is true since the weights of the edges are given by Lemma 6.4 as sums of weights of ends in a connected component of $\Gamma \setminus \{e\}$, and this differs for different edges. Let $\bar{e}$ denote the edge of weight $\pm\delta$. Naturally, $\bar{e}$ has weight $\delta$ for graphs in $\mathcal{T}^+(C_1)$ and $-\delta$ for graphs in $\mathcal{T}^+(C_2)$.

Therefore, we finally have:

$$\begin{aligned}
\mathrm{WC}_\delta = P_1 - P_2 &= \sum_{\Gamma \in (\mathcal{T}^+(C_1) \setminus \mathcal{T}^+(C_2))} \delta \prod_{e' \neq \bar{e}} \omega(e') \\
&\quad - \sum_{\Gamma \in (\mathcal{T}^+(C_2) \setminus \mathcal{T}^+(C_1))} (-\delta) \prod_{e' \neq \bar{e}} \omega(e') \\
&= \sum_{\Gamma \in \mathcal{T}^\delta} \delta \cdot \prod_{e' \neq \bar{e}} \omega(e').
\end{aligned}$$

□

We use the following notation: $\mu_I := (\mu_i)_{i \in I}$ is the subpartition of $\mu$ containing only the $\mu_i$ with index $i \in I$ (and $\nu_J$ analogously). Let $r = k + l - 2$, $r_1 = \#I + \#J - 1$ and $r_2 = k - \#I + l - \#J - 1 = r - r_1$.

**Theorem 6.10** [13] *For the wall $\delta = 0$ and the chambers $C_1$ and $C_2$ the wall crossing equals*

$$\mathrm{WC}_\delta(\mu,\nu) = \binom{r}{r_1, r_2} \cdot \delta \cdot H^0\big(\mu_I, (\nu_J, \delta)\big) \cdot H^0\big((\mu_{I^c}, \delta), \nu_{J^c}\big). \tag{7}$$

The idea of the proof is simple: cutting a graph contributing to the wall crossing along the special edge $\bar{e}$ disconnects the graph into two components, which we can interpret as graphs contributing to the Hurwitz numbers on the right hand side of the formula. Conversely, given any pair of graphs contributing to the Hurwitz numbers on the right hand side, we can glue them to create a graph contributing to the wall crossing. The binomial coefficient gives the number of ways to merge the total orderings of the vertices of the two components into a total order of the glued graph. We now make this precise. We maintain all notation from Lemma 6.9.

*Proof* The inequalities defining $C_1$ induce inequalities for the entries of the partitions $(\mu_I, (\nu_J, \delta))$ resp. $((\mu_{I^c}, \delta), \nu_{J^c})$. These inequalities determine chambers $C_1^1$ (resp. $C_1^2$) in the vector spaces $\mathbb{R}^{\#I+\#J+1}$ (resp. $\mathbb{R}^{k-\#i+l-\#J+1}$) with coordinates $(\mu_I, (\nu_J, \delta))$ resp. $((\mu_{I^c}, \delta), \nu_{J^c})$.

A graph $\Gamma_1 \in \mathcal{T}^+(C_1^1)$ has $r_1$ vertices. Any $\Gamma_2 \in \mathcal{T}^+(C_1^2)$ has $r_2$ vertices. Since $\Gamma_1 \cup \Gamma_2$ is disconnected, there are $\binom{r}{r_1,r_2}$ orderings of the $r = r_1 + r_2$ vertices of $\Gamma_1 \cup \Gamma_2$ that refine the orderings of $\Gamma_1$ and $\Gamma_2$ and are compatible with the directions of edges.

Let $\tilde{\mathcal{T}}$ denote the set of triples $(\Gamma_1, \Gamma_2, \tilde{O})$ with $\Gamma_1 \in \mathcal{T}^+(C_1^1)$, $\Gamma_2 \in \mathcal{T}^+(C_1^2)$ and $\tilde{O}$ is a choice of an ordering of the $r$ vertices of $\Gamma_1 \cup \Gamma_2$ as in the previous paragraph. Then:

$$\begin{aligned}
&\binom{r}{r_1, r_2} \cdot \delta \cdot H^0\big(\mu_I, (\nu_J, \delta)\big) \cdot H^0\big((\mu_{I^c}, \delta), \nu_{J^c}\big)\\
&\quad = \binom{r}{r_1, r_2} \cdot \delta \cdot \Bigg(\sum_{\Gamma_1 \in \mathcal{T}^+(C_1^1)} \varphi_{\Gamma_1}\Bigg) \cdot \Bigg(\sum_{\Gamma_2 \in \mathcal{T}^+(C_1^2)} \varphi_{\Gamma_2}\Bigg)\\
&\quad = \binom{r}{r_1, r_2} \cdot \delta \cdot \sum_{(\Gamma_1,\Gamma_2) \in \mathcal{T}^+(C_1^1) \times \mathcal{T}^+(C_1^2)} \varphi_{\Gamma_1} \cdot \varphi_{\Gamma_2}\\
&\quad = \delta \cdot \sum_{(\Gamma_1,\Gamma_2,\tilde{O}) \in \tilde{\mathcal{T}}} \varphi_{\Gamma_1} \cdot \varphi_{\Gamma_2}.
\end{aligned}$$

The theorem now follows from the claims:

**Claim 1:** There are natural bijections

$$\mathrm{Cut} : \mathcal{T}^\delta \leftrightarrow \tilde{\mathcal{T}} : \mathrm{Glue}.$$

**Claim 2:** The polynomial contributions of $\Gamma$ and $\mathrm{Cut}(\Gamma)$ to the two sides of (7) coincide.

Let $\Gamma \in \mathcal{T}^\delta$. If the graph $\Gamma$ is in $\mathcal{T}^+(C_1)$, then cut the special edge $\bar{e}$ to produce two graphs $\Gamma_1$ and $\Gamma_2$. If $\Gamma \in \mathcal{T}^+(C_2)$, then the new ends obtained by cutting $\bar{e}$ face the wrong way. The first graph would contribute to $H^0((\mu_I, \delta), \nu_J)$ instead of $H^0(\mu_I, (\nu_J, \delta))$, as desired. Hence we flip the direction of these ends. We now argue that this procedure does not produce any sources or sinks.

Let $v$ be the vertex at the source of $\bar{e}$, and suppose that both other edges incident to $v$ are oriented with $v$ as their target, so that reversing $\bar{e}$ would produces a sink. Recall that $\bar{e}$ is labeled with $\delta$; let $\alpha$ and $\beta$ be the labelings of the other two edges. The direction of the edges and the balancing condition then implies that $\alpha + \beta = \delta$. The hyperplanes defining chamber $C_1$ contain among them $\alpha > 0$ and $\beta > 0$, which together imply that $\delta > 0$, and so we see that the chamber $C_1$ does not border on the wall $\delta$ in codimension 1—to cross the wall $\delta = 0$ we must first cross one of the walls $\alpha = 0$ or $\beta = 0$. This is a contradiction, and so we see that one of the other edges incident to $v$ must have had $v$ as its source, and so flipping $\bar{e}$ does not create a sink. The argument for a source is analogous.

The data $\Gamma_1, \Gamma_2$, together with the ordering $\tilde{O}$ of all $r$ vertices that is naturally inherited by the ordering of the vertices in $\Gamma$, defines $\mathrm{Cut}(\Gamma)$.

We now define *Glue*. Let $(\Gamma_1, \Gamma_2, \tilde{O}) \in \tilde{\mathcal{T}}$. Glue $\Gamma_1$ and $\Gamma_2$ along the two ends with weight $\delta$, producing a new inner edge $\bar{e}$. The ordering $\tilde{O}$ of the $r$ vertices in particular orders the two vertices adjacent to $\bar{e}$ and thus determines the direction of $\bar{e}$. If the direction of $\bar{e}$ agrees with the original directions of the two ends we obtain a graph in $\mathcal{T}^+(C_1)$, otherwise the glued graph is in $\mathcal{T}^+(C_2)$. In the second case we must make sure no sources or sinks are created in changing the direction of $\bar{e}$. However, by definition $\Gamma_1$ belongs to chamber $C_1^1$, and so we can view it as a portion of a graph from chamber $C_1$ and follow the same argument we used in defining the map *Cut*, and similarly $\Gamma_2$ is a graph contributing to chamber $C_1^2$.

Thus, we see the glued graph together with the ordering $\tilde{O}$ of all $r$ vertices belongs to $\mathcal{T}^\delta$, and is defined to be $\mathrm{Glue}(\Gamma_1, \Gamma_2, \tilde{O})$.

The two maps *Cut* and *Glue* are inverses of each other and thus give bijections of the two sets, as stated in Claim 1. Furthermore, $\Gamma$ contributes $\delta \cdot \prod_{e' \neq \bar{e}} \omega(e')$ to the left hand side of (7). The image $\mathrm{Cut}(\Gamma)$ contributes

$$\begin{aligned} \delta \cdot \varphi_{\Gamma_1} \cdot \varphi_{\Gamma_2} &= \delta \cdot \prod_{e' \text{ int. edge of } \Gamma_1} \omega(e') \cdot \prod_{e' \text{ int. edge of } \Gamma_2} \omega(e') \\ &= \delta \cdot \prod_{e' \neq \bar{e},\ e' \text{ int. edge of } \Gamma} \omega(e') \end{aligned}$$

to the right hand side of (7), proving Claim 2 and Theorem 6.10. □